# A Methodology for Projection-Based Model Reduction with Black-Box High-Fidelity Models

S. Ashwin Renganathan[a], Yingjie Liu[b], and Dimitri N. Mavris[c]
*Georgia Institute of Technology, 270 Ferst Dr, Atlanta, GA 30313*

**This paper presents a methodology that enables projection-based model reduction for black-box high-fidelity models such as commercial CFD codes. The methodology specifically addresses the situation where the high-fidelity model may be a black-box but there is complete knowledge of the governing equations. The main idea is that the linear operator matrix, resulting from the discretization of the linear differential terms is approximated directly using a suitable discretization method such as the Finite Volume Method and requires only the computational grid as input. In this regard, the governing equations are first cast in terms of a set of scalar observables of the state variables, leading to a linear set of equations. By applying the snapshots of the observables to the discrete linear operator, a right-hand side vector is obtained, providing the necessary system matrices for the projection step. This way an offline database of ROMs is generated for various parameter snapshots which are then interpolated online to predict the state for new parameter instances. Finally, the reduced order model is posed as a non-linear constrained optimization problem that can be solved at a significantly cheaper cost compared to the full order model. The method is successfully demonstrated on a canonical non-linear parametric PDE with exponential non-linearity, followed by the compressible inviscid flow past the NACA0012 airfoil. As a first step, this paper focuses only on establishing feasibility of the method.**

[a] Ph.D. Candidate, School of Aerospace Engineering, Student Member. email: ashwinsr@gatech.edu
[b] Professor, School of Mathematics, Non-Member.
[c] Regents Professor, School of Aerospace Engineering, AIAA Fellow

## I. Introduction

Aerospace design optimization typically involves a high dimensional design space with multiple local optima. To ensure that the globally optimum design is obtained, in addition to having robust and accurate optimization methods, it is also important to use accurate high-fidelity models to search the design space. The search for global optimum using high-fidelity models in engineering design optimization is typically done using Surrogate Based Optimization (SBO) [1–3]. Surrogate models or Metamodels are developed by regressing a training data set generated via high fidelity models at a pre-determined set of points and the optimum is searched within the surrogate model. There is some sophistication available to this approach where the surrogate model can be adaptively developed such that more samples are concentrated in regions where the optimum is likely to be found [1]. Such models are generally developed using data-fit approaches such as Response Surface Methodology (RSM), Artificial Neural Networks (ANN) or interpolating models such as Kriging and Radial Basis Functions. While they are easy to implement and assume simple functional forms that are computationally cheap to evaluate, the limitations are that such models (i) suffer from the curse of dimensionality and (ii) do not account for the physics of the problem and hence are not guaranteed to provide a certain definite level of accuracy. Another class of surrogate models that differ from the data-fit models by being customizable to the physics of the high-fidelity models are called *physics-based metamodels*. Their name is due to the fact that their formulation ensures they satisfy the governing equations of the system [4].

Reduced Order Models (ROMs) are a class of physics-based metamodels that develop a low-dimensional representation of the high-fidelity model, the Full Order Model (FOM). Therefore, the spatial and temporal variation of the state variables of the system can be obtained by solving a much smaller version of the FOM which is computationally cheaper, while still satisfying the governing equations. Model reduction for non-linear parametric systems is specifically suited for design optimization and has been applied to a variety of aerospace design problems in the past decade including aerodynamic inverse design and missing data re-construction [5, 6], probabilistic aerodynamic analysis [7], aero-elastic applications [7–16], aero-thermo-elastic applications [17, 18], reacting flow modeling [19], hypersonic thermal protection system design [20, 21], turbomachinery [22–24], rotary wing and prop-rotor aerodynamics [25–27] and supersonic flow modeling [28–31]. The current state of the art of ROM development is *intrusive* in nature since it requires access to the governing equations of the FOM to construct the ROM. However, when the FOM is solved using a black-box code such as commercial Computational Fluid Dynamics (CFD) packages, the governing equations are inaccessible, posing a hurdle to the ROM development. This paper aims to address this limitation.

The fundamental assumption behind ROMs is that most of the variance of a high-dimensional flow field is explained by a low-dimensional vector space spanned by a set of orthonormal basis vectors. Therefore, by projecting the original governing equations onto the space spanned by the basis set results in dimensionality reduction in terms of the number of unknowns being solved for. The process of extracting the basis set is achieved via methods typically based on the Singular Value Decomposition (SVD) or Krylov Subspace Methods [32]. Proper Orthogonal Decomposition (POD) [33, 34] is a popular SVD-based method that are particularly suitable for non-linear model reduction and is briefly discussed in the following section.

### A. Proper Orthogonal Decomposition

POD was originally introduced in the context of turbulent flow modeling by Holmes et al [34], where it was used to characterize the coherent structures in the flow from wind tunnel measurements. POD has a special characteristic of *optimality* in that it provides the most efficient means to capture the dominant components of a process [35]. Given a state variable $\mathbf{u} \in \mathbb{R}^N$ which may be the numerical solution of a PDE on a computational mesh of size $N$, the POD expresses $\mathbf{u}$ as the linear combination of a finite number of $k$ orthonormal basis vectors $\phi_i \in \mathbb{R}^N$. That is,

$$\mathbf{u} \approx \sum_{i=1}^{k} \tilde{\mathbf{u}}(i)\phi_i \tag{1}$$

where, $\tilde{\mathbf{u}}(i)$ is the $i^{th}$ component of $\tilde{\mathbf{u}} \in \mathbb{R}^k$ and are the coefficients of the basis expansion. Denoting $\mathbf{\Phi}_k = [\phi_1, ..., \phi_k]$, the Equation 1 can be written as

$$\mathbf{u} \approx \mathbf{\Phi}_k \tilde{\mathbf{u}} \tag{2}$$

Therefore, by substituting Equation 2 onto the original governing equations, we may solve for $\tilde{\mathbf{u}}$ instead of $\mathbf{u}$, which is more efficient since $k << N$. The POD basis is determined such that it minimizes the error between the state variable and its orthogonal projection onto $\Phi_k$. Given $M$ snapshots of the state variables $\mathbf{U} = [\mathbf{u}^1, ..., \mathbf{u}^M] \in \mathbb{R}^{N\times M}$ which may be obtained by solving the FOM for various independent parameter combinations, the POD basis satisfies

$$\phi_i := \arg\min \sum_{j=1}^{M} \|\mathbf{u}^j - \sum_{i=1}^{k} ({\mathbf{u}^j}^T \phi_i)\phi_i\|_2^2 = \arg\min \sum_{j=1}^{M} \|\mathbf{u}^j - (\mathbf{\Phi}_k \mathbf{\Phi}_k^T)\mathbf{u}^j\|_2^2 \tag{3}$$

where $(\mathbf{\Phi}_k \mathbf{\Phi}_k^T)\mathbf{u}^j$ is the orthogonal projection of $\mathbf{u}^j$ onto $\mathbf{\Phi}_k$ and $\mathbf{\Phi}_k^T \mathbf{\Phi}_k = \mathbf{I}$. It can be shown that [35, 36] the solution to the above equation is the same as the left singular vectors of the snapshot matrix. That is,

$$\mathbf{U} = \mathbf{V}\mathbf{\Sigma}\mathbf{W}^T \tag{4}$$

then $\mathbf{\Phi}_k$ represents the first $k$ columns of $\mathbf{V} \in \mathbb{R}^{N\times N}$. The $L_2$ error in approximation of the state variables due to the POD basis expansion is then given as

$$\sum_{j=1}^{M} \left\|\mathbf{u}^j - (\mathbf{\Phi}_k \mathbf{\Phi}_k^T)\mathbf{u}^j\right\|_2^2 = \sum_{i=k+1}^{M} \sigma_i^2 \tag{5}$$

where $\sigma_i$ are the singular values of the $i^{th}$ column of $\mathbf{V}$ and are also the $i^{th}$ diagonal element of $\mathbf{\Sigma}$. The POD step is followed by the projection step where the FOM governing equations are projected onto the reduced set of POD basis thereby reducing the original $N \times N$ system into a $k \times k$ system where $k$ is the dimension of the reduced POD basis set. This is briefly discussed in the following section.

### B. Projection Step

The projection step projects the original FOM onto the low-dimensional subspace spanned by the POD basis vectors, forming the ROM [35]. Consider the discrete representation of a steady, non-linear, parametric system that is the result of the discretization of a PDE

$$\mathbf{A}(\theta)\mathbf{u} = \mathbf{f}(\mathbf{u}) \tag{6}$$

where $\mathbf{A}(\theta) \in \mathbb{R}^{N \times N}$ is the linear differential operator that arises due to the discretization of linear terms, $\theta \in \mathbb{R}^p$ is a vector of design parameters and $\mathbf{f}(\mathbf{u}) \in \mathbb{R}^{N \times 1}$ is the non-linear operator that arises due to the discretization of the non-linear terms and also lumps the boundary condition discretization terms and source terms if present. This represents the full-order system with $N$ unknowns. The projection step begins by realizing that the residual of the FOM is orthogonal to an appropriately chosen test basis, $\mathbf{\Psi_k}$. Removing the $\theta$ for convenience of notation, this is equivalent to

$$\mathbf{\Psi}_k^T(\mathbf{A}\mathbf{u} - \mathbf{f}(\mathbf{u})) = 0 \tag{7}$$

Since the reduced state variable $\tilde{\mathbf{u}} \approx \mathbf{\Phi}_k^T\mathbf{u}$, the above equation can be written as

$$\mathbf{\Psi}_k^T\mathbf{A}\mathbf{\Phi}_k\tilde{\mathbf{u}} = \mathbf{\Psi}_k^T\mathbf{f}(\mathbf{\Phi}_k\tilde{\mathbf{u}}) \tag{8}$$

defining the reduced matrix $\tilde{\mathbf{A}} = \mathbf{\Psi}_k^T\mathbf{A}\mathbf{\Phi}_k \in \mathbb{R}^{k \times k}$,

$$\tilde{\mathbf{A}}\tilde{\mathbf{u}} = \mathbf{\Psi}_k^T\mathbf{f}(\mathbf{\Phi}_k\tilde{\mathbf{u}}) \tag{9}$$

The above equation represents a reduced system with $k << N$ unknowns that can be solved efficiently. The computation of the non-linear term still involves operations in $\mathcal{O}(N)$ and renders the computation inefficient since it has to be evaluated repeatedly in an iterative procedure such as the Newton's method. However, this can be overcome with the Discrete Empirical Interpolation Method (DEIM) due to Chaturantabut & Sorensen [37] to compute it efficiently and is adopted in this paper.

The main motivation for this work comes from the fact that the projection step explained above requires the system matrices $(\mathbf{A}, \mathbf{f})$ which are not available in situations where the full model is a black-box, such as commercial CFD codes. Therefore, projection-based ROMs are not feasible in such situations. There have been other attempts to address or work around this situation in the literature. One approach is to replace the projection step with an interpolation step, where the reduced state variables are directly interpolated for time and/or parametric changes using a suitable interpolation method. See [38, 39] for instance, where $2^{nd}$ order Taylor series expansion, Smolyak sparse grid collocation and RBFs are used for interpolation. On the other hand, data-driven approaches exist to enable projection-based model reduction for non-linear dynamic systems. A linear operator that maps the time-evolution of the state variables using trajectory data can be approximated, which can further be used for model reduction. Such an approach derives a finite-dimensional approximation of the so-called *Koopman Operator* [40, 41]. A similar approach has been taken by [42] to infer the reduced operator matrices directly, which was demonstrated on linear PDEs and PDE's with low-order polynomial non-linearity. In this paper we aim to address situations where projection-based ROMs are of interest while there is no trajectory data available and the system is parametric.

This work proposes a method that enables projection-based construction of ROMs with black-box full models for steady, non-linear parametric systems with generalized non-linearities. It specifically addresses the situation where there is knowledge of the governing equations, however there is no access to the source code which is typically required to perform model reduction. It first represents the governing equation in terms of a set of *observables*, drawing from the Koopman theory of partial differential equations [40]. Such a representation transforms the governing equation into a higher dimensional but linear system. Secondly, the linear operator matrix $\mathbf{A}$ is approximated via a direct discretization of the transformed equation. The method recognizes that $\mathbf{A}$ is the discrete version of a few standard linear differential operators such as the gradient, divergence, laplacian and curl. Additionally, $\mathbf{A}$ requires only the computational grid and the parameters of the system which are always known. For a given snapshot of the FOM (in terms of the observables) and the matrix $\mathbf{A}$, the right-hand side, $\mathbf{f}$ of the governing equations (which lumps the boundary condition and source terms) is recovered. The overall method is summarized in the Figure 1. Note that the $\mathbf{f}$ now is independent of the state since all the non-linear terms are cast in terms of an observable and hence is suitable for a direct element-wise interpolation for parameter changes. This allows one to generate the system matrices $\mathbf{A}$ and $\mathbf{f}$ that are necessary for projection-based reduced order modeling.

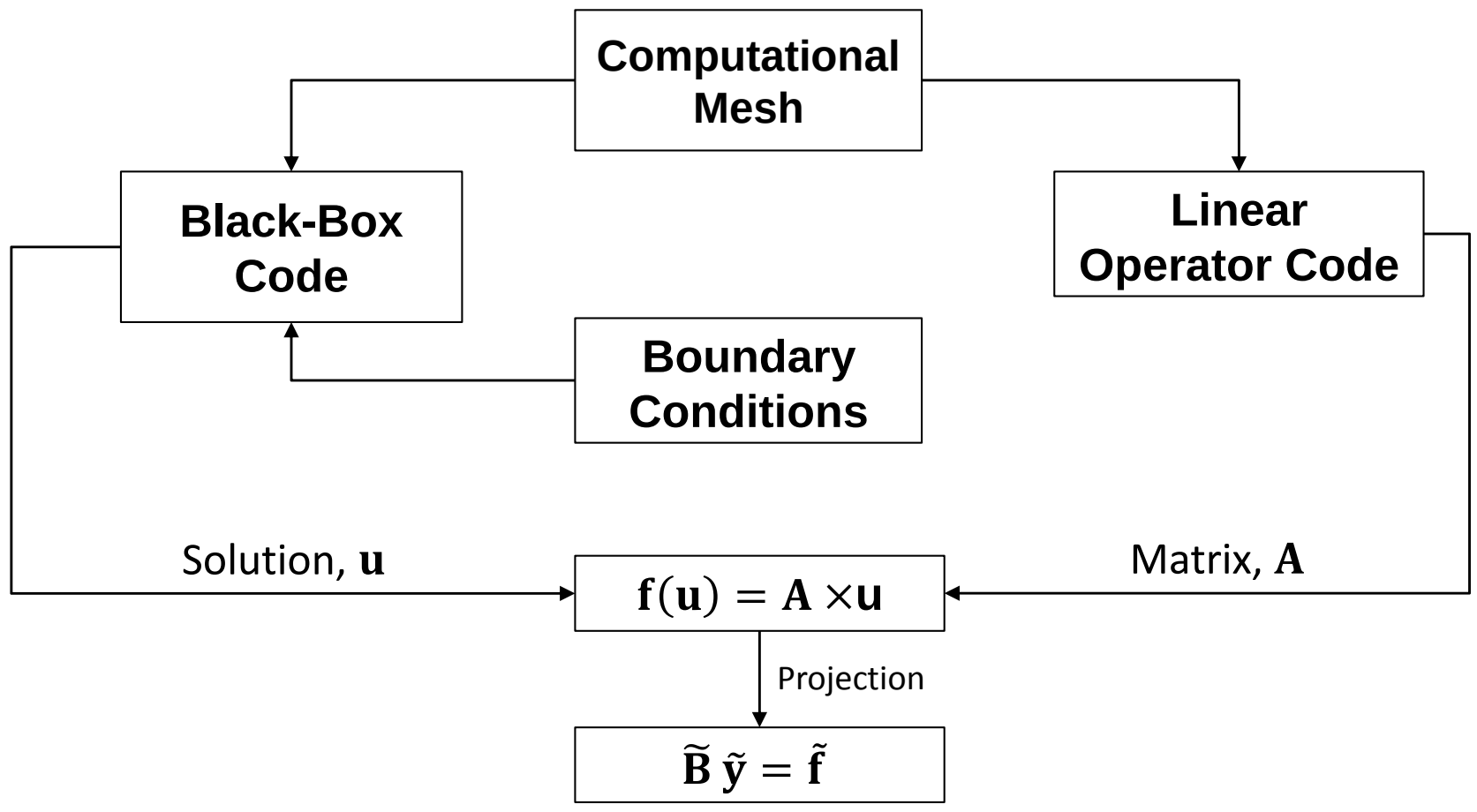

**Figure 1 Flowchart of overall workflow**

The rest of this paper is organized as follows. The next section describes the methodology in detail, following which the computational complexity is discussed and an overall algorithm is provided. Following that the results of the test problems used to demonstrate the method are provided followed by a discussion. The paper concludes with a summary of the results and an outlook on future work.

## II. Methodology

The present approach first begins by representing the governing equations in terms of a set of scalar observables. As mentioned before, this step draws from the Koopman Theory which is briefly described in the Appendix A and in greater detail in [40, 43, 44]. The main purpose of this step is that the resulting equations lead to a form $\mathbf{A}\mathbf{y} = \mathbf{f}$ where both $\mathbf{A}$ and $\mathbf{f}$ are independent of the state and hence are amenable to an element-wise interpolation to adapt for parametric changes. Following this, the finite volume method is described which is used to obtain the linear operator $\mathbf{A}$ by directly discretizing the linear terms of the transformed equation. This section details the steps involved in the proposed methodology.

## A. Representation in observables space

Consider a steady non-linear system of the form

$$\mathcal{N}(\mathbf{u}) = 0 \tag{10}$$

where $\mathcal{N}$ represents a non-linear operator on the state variable $\mathbf{u}$. Let $g(\mathbf{u})$ represent an observable that is a function of the state variable, $\mathbf{u}$. We state that

$$\mathcal{N}(\mathbf{u}) \equiv \mathcal{L}[g(\mathbf{u})] \tag{11}$$

where, $\mathcal{L}$ is a linear operator acting on the observables. Discretizing the above equation, we get

$$\mathcal{L}[g(\mathbf{u})] \approx \mathbf{A}g(\mathbf{u}) + \mathbf{b}_a = 0 \tag{12}$$

Where, $\mathbf{b_a}$ represents a vector that arises due to the discretization of the boundary conditions, lumps the source terms if present and is also the RHS of the FOM. Note again that since each non-linear term is transformed into an observable, $\mathbf{f}$ (in $\mathbf{Au} = \mathbf{f}$) no longer depends on the state $\mathbf{u}$ and is equal to $-\mathbf{b_a}$ in Equation 12. Also, we intend to develop the ROM for Equation 12 and hence the snapshots are collected in the observable space, $g(\mathbf{u})$. Therefore by discretizing the linear terms of the governing equations we obtain $\mathbf{A}$ which is then applied onto the snapshots to obtain the RHS, $\mathbf{f}$.

## B. Model Reduction

We begin by representing the FOM as

$$\mathbf{Ay} = \mathbf{f} \tag{13}$$

where, $g(\mathbf{u})$ is replaced with $\mathbf{y}$ and $\mathbf{A}$, $\mathbf{f}$ can have parametric dependence but are independent of the state. A side effect of representing the FOM in terms of observables is that it can increase the dimensionality of the system. For instance consider a nonlinear system with 1 state variable $\mathbf{u}$ and 2 observables $\mathbf{y} = [g_1(\mathbf{u})\ g_2(\mathbf{u})]^T$. Now, the system can be written as

$$\begin{bmatrix} \mathbf{A_1} & \mathbf{A_2} \end{bmatrix} \begin{bmatrix} g_1(\mathbf{u}) \\ \\ g_2(\mathbf{u}) \end{bmatrix} = \mathbf{f} \tag{14}$$

where $\mathbf{A} = [\mathbf{A_1}\ \mathbf{A_2}]$. Since each of $g_1(\mathbf{u})$ and $g_2(\mathbf{u})$ have the same dimensions of the state variable, $\mathbf{u}$, this leads to an under-determined system. Therefore, in order to ensure a unique solution, additional constraints are required to be added to the solution of the ROM, which are discussed later on in this section.

Denoting $\phi_1 \in \mathbb{R}^{k_1}$ and $\phi_2 \in \mathbb{R}^{k_2}$ to be the reduced POD bases extracted from snapshots of $g_1$ and $g_2$ respectively, the trial basis matrix is defined as

$$\Phi = \begin{bmatrix} \phi_1 & \\ & \phi_2 \end{bmatrix} \in \mathbb{R}^{2N \times k} \tag{15}$$

where $k = k_1 + k_2$ and $\tilde{\mathbf{y}} = \Phi_k^T \mathbf{y}$. Since the full system is non-square (due to introduction of observables), a suitable choice for the test basis for projection is $\mathbf{\Psi}_k = \mathbf{A}\mathbf{\Phi}_k$. Note that this choice of the test basis is equivalent to a galerkin projection ($\Psi_k = \Phi_k$) on the normal equations. i.e. on $\mathbf{A}^T\mathbf{A}\mathbf{y} = \mathbf{A}^T\mathbf{f}$. Let $\mathbf{B} = \mathbf{A}^T\mathbf{A}$ , then the projection leads to

$$\Phi_k^T \mathbf{B} \Phi_k \tilde{\mathbf{y}} = \Phi_k^T \mathbf{A}^T \mathbf{f} \tag{16}$$

Setting $\tilde{\mathbf{f}} = \Phi_k^T \mathbf{A}^T \mathbf{f} \in \mathbb{R}^k$ and $\tilde{\mathbf{B}} = \Phi_k^T \mathbf{B} \Phi_k \in \mathbb{R}^{k \times k}$, this leads to the reduced order model

$$\tilde{\mathbf{B}} \tilde{\mathbf{y}} = \tilde{\mathbf{f}} \tag{17}$$

The same derivation can be extended for arbitrary number of observables without modification. The ROM given by Equation 17 represents a $k \times k$ system which is rank-deficient since it was obtained through an outer product of two rectangular matrices. Therefore, in order to enforce the equations of states and hence uniqueness of the solution, the ROM is posed as a constrained optimization problem as shown below

$$\begin{aligned} \underset{\tilde{\mathbf{y}}}{\text{minimize}} \quad & \frac{1}{2}\|\tilde{\mathbf{B}}\tilde{\mathbf{y}} - \tilde{\mathbf{f}}\|_2^2 \\ \text{s.t.} \quad & h(\mathbf{y}) = 0 \end{aligned} \tag{18}$$

where $h(\mathbf{y})$ is a non-linear function that represents the relationship between $\mathbf{y}_1$ and $\mathbf{y}_2$ which depends on the choice of $g_1$ and $g_2$ and hence is problem dependent. The main hypothesis of this work is that the ROM given by Equation 18 still approximately satisfies the governing equations and this is verified in the section IV. The optimization problem in Equation 18 needs special treatment to handle the non-linear constraint which still depends on the full state of observables, and is efficiently done using the DEIM which is briefly reviewed in Appendix B but further details can be obtained from [37].

### C. Finite Volume Discretization

As previously mentioned, the matrix $\mathbf{A}$ in the present method is directly approximated via discretization of the linear terms. A finite volume based approach [45] is used for the discretization due to its suitability to complex geometries with arbitrarily shaped cells. The final matrix $\mathbf{A}$ itself is composed of several matrices each of which represents the discretized form of the linear differential terms present in the governing equations of the FOM.

In this work, we take the cell-centered finite volume approach where the dependent variable is stored in the cell centers. In this subsection as an illustration we briefly review the discretization of the 2D diffusion operator $\nabla \cdot (\Gamma\nabla)$ which reduces to the laplacian, $\nabla^2$ when $\Gamma = 1$.

Consider the linear diffusion term is given by

$$\mathcal{L}u = \nabla.(\Gamma\nabla u) \tag{19}$$

where the diffusion coefficient $\Gamma = \Gamma(x, y)$ is independent of $u$. Integrating the above term over a cell volume and applying the Gauss-Divergence theorem we get,

$$\int \nabla.(\Gamma\nabla u)dV = \oint (\Gamma\nabla u).\hat{n}ds \tag{20}$$

where $ds$ is the face area and $\hat{n}$ is the local surface normal of the face. Assuming that the computational mesh consists of only polygonal faces, each face of a cell has a unique surface normal and hence the above integral can be written as

$$\oint (\Gamma \nabla u).\hat{n} ds = \sum_f \Gamma_f (\nabla u)_f.\hat{n}_f A_f \tag{21}$$

where the summation is over the faces of a given cell (subscripted by $f$), $\Gamma_f$ is the diffusion coefficient at the interface $f$ and $A_f$ is the face area (length in 2D). The face-center values of diffusion coefficient are obtained from the cell-center values via a linear interpolation (where the interpolants $w_f$ are inverse distance weighted) as

$$\Gamma_f = w_f \Gamma_0 + (1 - w_f)\Gamma_1$$

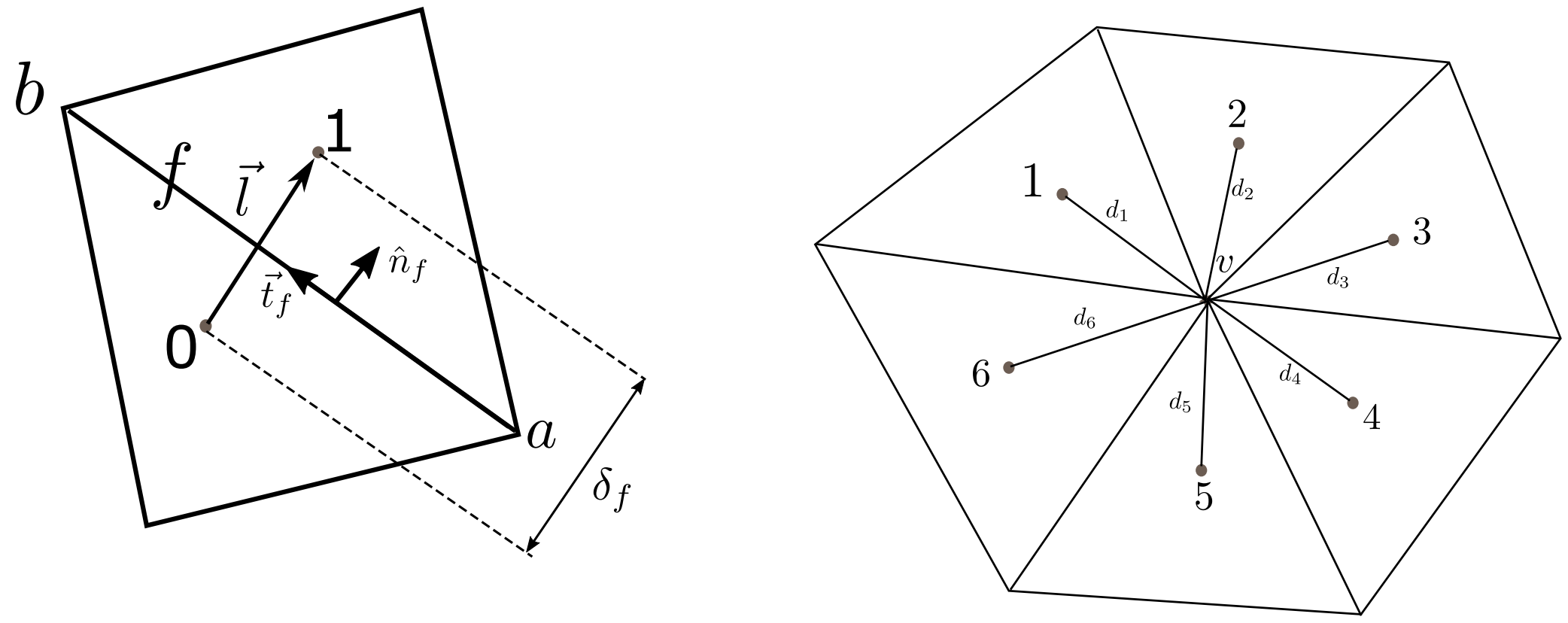

(a) Computation of fluxes at cell interface

(b) Distance-weighted interpolation of node values

**Figure 2 Finite volume cells**

Now it is a matter of discretizing the RHS of Equation 21. Consider two adjacent cells as shown in Figure 2a. Here, the center of the cell upon which the discretization is carried out is denoted as '0' and that of its neighboring cell as '1'. The cell centered values of $u$ are denoted as $u_0$ and $u_1$ for each of these cells respectively. The vector connecting the cell centers is denoted $\vec{l}$ and that connecting the vertices of the interface $f$ is $\vec{t}_f$ and $\hat{n}_f$ is the unit surface normal. Note that the $(\vec{\ })$ notation implies vector and the $\hat{( )}$ notation implies a unit vector. Finally, $\delta = \vec{l}.\hat{n}_f$ and since $\hat{n}$ and $\hat{t}$ form an orthogonal set of coordinate vectors, the gradient can be written as

$$(\nabla u) = [\nabla u.\hat{n}]\hat{n} + [\nabla u.\hat{t}]\hat{t} \tag{22}$$

At the face $f$, taking the dot product of $(\nabla u)_f$ with $\vec{l}$ we get

$$(\nabla u)_f.\vec{l} = [\nabla u.\hat{n}_f]\delta + [\nabla u.\hat{t}_f]\hat{t}_f.\vec{l} \tag{23}$$

Using taylor series to expand $u_0$ and $u_1$ about $u_f$ and subtracting we get

$$u_1 - u_0 \approx \left(\frac{\partial u}{\partial x}\right)_f (x_1 - x_0) + \left(\frac{\partial u}{\partial y}\right)_f (y_1 - y_0) = (\nabla u)_f.\vec{l} \tag{24}$$

i.e.

$$(\nabla u)_f.\vec{l} \approx u_1 - u_0 \tag{25}$$

In the above equation, the coordinates of the cell centers $'0'$ and $'1'$ are $(x_0, y_0)$ and $(x_1, y_1)$ respectively. Therefore the term $(\nabla u)_f.\hat{n}_f$ is given as

$$(\nabla u)_f.\hat{n}_f = \frac{u_1 - u_0}{\delta_f} - \frac{[\nabla u)_f.\hat{t}_f]\hat{t}_f.\vec{l}}{\delta_f} \tag{26}$$

In the Equation 26 the second term is called the 'tangential flux' and denoted as $J_T$ where $\frac{J_T}{\delta_f} = \frac{[(\nabla u)_f.\hat{t}_f]\hat{t}_f\vec{l}}{\delta_f}$. For structured cartesian meshes certain unstructured meshes such as those with equilateral triangular cells, the tangential flux term vanishes and the above expression reduces to standard central differencing. So now we proceed to discretize $J_T$. We begin by realizing that the Equation 25 can be re-written as

$$(\nabla u)_f.\hat{l} \approx \frac{u_1 - u_0}{|\vec{l}|} \tag{27}$$

where $|\vec{l}|$ is the inter-cell distance. Similarly, we write

$$(\nabla u)_f.\hat{t}_f \approx \frac{u_a - u_b}{|\vec{t}_f|} \tag{28}$$

$|\vec{t}_f|$ is the length of the face connecting vertices $a$ and $b$ in Figure 2a. Therefore, the tangential flux $J_T$ can be written as

$$J_T = \left[\frac{u_a - u_b}{|\vec{t}_f|}\right] \hat{t}_f.\vec{l} \tag{29}$$

In the above equation, since we know the vectors $\hat{t}_f$ and $\vec{l}$, their dot product is directly obtained. Thus putting it all together into Equation 21,

$$\oint (\Gamma \nabla u).\hat{n} ds = \sum_f \Gamma_f \left[\frac{u_{k(f)} - u_0}{\delta_f} - \left[\frac{u_a - u_b}{\delta_f |\vec{t}_f|}\right] \hat{t}_f.\vec{l}\right] A_f \tag{30}$$

where, $u_{k(f)}$ represents cell centered values of all the neighbors of cell '0'. The terms $(u_a, u_b)$ in the evaluation of the tangential flux are node-based values (at nodes $a$ and $b$) and are explicitly treated as the distance-weighted average of all the neighboring cell-centered values, as shown in Figure 2b. i.e. $u$ at a vertex $v$ is given as

$$u_v = \sum_i w_{v,i} u_i \tag{31}$$

where $i = 1, 2, 3, ...$ represent the cells surrounding vertex $v$. The actual number of neighboring cells for a vertex depends on the type of mesh and its location near or away from boundaries. The interpolant $w_{v,i}$ are given by

$$w_{v,i} = \frac{1/d_i}{\sum_i 1/d_i} \tag{32}$$

where $d_i$ is the distance of cell center of neighbor cell $i$ to node $v$.

It should be noted that the faces containing boundary condition information are assumed to be lumped by the RHS vector $\mathbf{f}$. In the finite-volume method, this approach is consistent for a pure

Neumann type boundary condition where the flux at the boundary is a known quantity but for other types of boundary conditions this introduces some error in the approximation of the matrix $\mathbf{A}$. However, with knowledge of the boundary conditions, a more accurate approximation of $\mathbf{A}$ can be obtained via the present approach.

### D. ROM Interpolation

The resulting matrices $\tilde{\mathbf{B}}$ and $\tilde{\mathbf{f}}$ are in general parameter dependent and they are unique for each snapshot. For new instances of the parameters, the matrices are interpolated as opposed to re-constructing them which is computationally expensive. As mentioned earlier, the matrices can be interpolated element-wise using methods in [46] or [47] for example. In order to focus on the overall model reduction method, this paper performs an element-wise interpolation in the original space using polynomials in the Lagrange form [48, 49]. This is briefly explained for the multivariate case in what follows.

We are interested in constructing a degree $n$ polynomial of an $m$-variate function $f(\mathbf{x}),\ \ \mathbf{x} \in \mathbb{R}^m$ which interpolates $\rho' = \binom{n+m}{m}$ points. Therefore, there needs to be atleast $\rho'$ snapshots to interpolate the reduced matrices. We are interested in polynomials of the form $\ell_i(\mathbf{x})$ where

$$f(\mathbf{x}) = \sum_{i=1}^{\rho'} f_i \ell_i(\mathbf{x}) \tag{33}$$

In Equation 33, $\ell_i(\mathbf{x_i}) = 1$ and $\ell_i(\mathbf{x} \neq \mathbf{x_i}) = 0$ which gives $f(\mathbf{x_i}) = f_i$ and hence interpolating the true function. The determination of $\ell(\mathbf{x})$ begins by representing $f(\mathbf{x})$ in the following coefficient form

$$f(\mathbf{x}) = \sum_{\mathbf{e}_i \cdot \mathbf{1} \leq n} \alpha'_{\mathbf{e}_i} \mathbf{x}^{\mathbf{e}_i} \tag{34}$$

where, $\mathbf{1} = [1, ..., 1]^T \ \in \mathbb{R}^m$. This leads to a system of $\rho'$ linear equations of the form $f_i = \sum_j \alpha'_{\mathbf{e}_j} \mathbf{x_i}^{\mathbf{e}_j}$ where $\mathbf{e_i} = (e_{1i}, \ldots, e_{mi})$ is a vector of non-negative integers such that the sum of the elements is $\leq n$, $\alpha'_{\mathbf{e_i}}$ are the coefficients and $\mathbf{x}^{\mathbf{e_i}} = \prod_{j=1}^{m} \mathbf{x_j}^{e_{ji}}$. Rather than determining a functional form of the interpolant $\ell_i(\mathbf{x})$, the coefficients may be determined as the solution of the linear system $\mathbf{M}\alpha' = f$ at the $\varrho'$ points where matrix $\mathbf{M}$ takes the form

$$\mathbf{M} = \begin{bmatrix} \mathbf{x_1}^{\mathbf{e_1}} & \dots & \mathbf{x_1}^{\mathbf{e}_{\varrho'}} \\ \vdots & \dots & \vdots \\ \mathbf{x}_{\varrho'}^{\mathbf{e_1}} & \dots & \mathbf{x}_{\varrho'}^{\mathbf{e}_{\varrho'}} \end{bmatrix} \tag{35}$$

In this work we use the $2^{nd}$ order, bi-variate ($m = 2,\ n = 2$) version of Equation 33 for interpolation and is illustrated as follows. Let $\theta_{train} \in \mathbb{R}^{M\times 2}$ be the parameters at which the $M$ snapshots are generated, $\hat{\theta} \in \mathbb{R}^2$ be the parameter value at which interpolation is performed and we are interested in interpolating the first element of the RHS vector $\mathbf{f}$ at $\hat{\theta}$, denoted by $\hat{\mathbf{f}}_1$. Then there are a total of $\varrho' = \binom{n+m}{n} = 6$ points are required to interpolate. This may be chosen as the $\varrho' < M$ nearest neighbors to $\hat{\theta}$ in $\theta_{train}$. Then,

$$\hat{\mathbf{f}}_1 = a^T \vartheta \tag{36}$$

where $\vartheta$ is the solution of $\mathbf{M}\vartheta = \mathbf{f}_1(i_{\varrho'})$ and $i_{\varrho'}$ are the indices of the $\varrho'$ points. $\mathbf{M} \in \mathbb{R}^{\varrho'\times\varrho'}$ is the matrix of basis functions evaluated at the $\varrho'$ points and is given by

$$\mathbf{M} = \begin{bmatrix} \vdots & \vdots & \vdots & \vdots & \vdots & \vdots \\ 1 & \theta_{train,1} & \theta_{train,2} & \theta_{train,1}\theta_{train,2} & \theta^2_{train,1} & \theta^2_{train,2} \\ \vdots & \vdots & \vdots & \vdots & \vdots & \vdots \end{bmatrix} \tag{37}$$

Finally, $a$ is a column vector of the basis functions evaluated at $\hat{\theta}$, i.e.

$$a = \begin{bmatrix} 1 & \hat{\theta}_1 & \hat{\theta}_2 & \hat{\theta}_1\hat{\theta}_2 & \hat{\theta}_1^2 & \hat{\theta}_2^2 \end{bmatrix}^T \tag{38}$$

The above illustration interpolates 1 element and is repeated for every element of a given vector or matrix.

### E. Overall Algorithm

The overall method is summarized here. Note that in the algorithm, $D$ represents the design space while $D_{train}$ represents a subset of $D$ used for snapshot generation (model training).

**Algorithm 1** Non-Intrusive Projection-Based Model Reduction

1: Choose $M$ snapshot locations for model training $\theta_i \in D_{train} \subset D$

2: Solve FOM and construct 'observables' snapshot matrix: $\mathbf{U} = [g(\mathbf{u}^1), ..., g(\mathbf{u}^M)],\ \mathbf{u}^i \in \mathbb{R}^N$

3: Extract the trial basis $\Phi_k$ via POD: $\mathbf{U} = \mathbf{V}\Sigma\mathbf{W}^T$, $\Phi_k = \mathbf{V}(:, 1:k)$

4: Construct system matrices

  **for** *every* $\theta_i$

   $\mathbf{A}_i \leftarrow$ Discretize Linear Operator

   $\mathbf{f}_i \leftarrow \mathbf{A}_i \times \mathbf{u}_i$

   $\tilde{\mathbf{B}}_i \leftarrow \Phi_k^T(\mathbf{A}_i^T\mathbf{A}_i)\Phi_k$

   $\tilde{\mathbf{f}}_i \leftarrow \Phi_k^T\mathbf{A}_i^T\mathbf{f}$

  **end for**

5: Prediction: for any $\theta' \notin D_{train},\ \theta' \in D$

  interpolate piece-wise $\tilde{\mathbf{B}}(\theta')$

  interpolate piece-wise $\tilde{\mathbf{f}}(\theta')$

6: Solve ROM:

$$\underset{\tilde{\mathbf{y}}}{\text{minimize}} \quad \frac{1}{2}\|\tilde{\mathbf{B}}\tilde{\mathbf{y}} - \tilde{\mathbf{f}}\|_2^2$$

$$\text{s.t.} \qquad h(\mathbf{y}) = 0$$

7: Project ROM onto FOM space: $\mathbf{y} = \Phi_k\tilde{\mathbf{y}}$

8: Map observables back to state variables

  $\mathbf{u} \leftarrow \mathbf{y}$

## III. Offline Computational Cost

The construction of the linear differential operator is dependent on knowing the connectivity between the various elements of the mesh. For instance, information about the cells that lie adjacent to a face, the faces that bound a cell, the vertices that make up each edge and the cells that share a vertex are examples of connectivity that encode an unstructured mesh that needs to be known to discretize a PDE on it. Surface normals are also required to compute the component of normal and tangential flux on each cell face. Additionally, the interpolants required to convert between cell-center, face-center and node values are also to be computed. The current section focuses only on estimating the computational complexity and storage for discretizing the linear operator matrix, given all the other aforementioned information since that dominates the overall offline cost of the proposed approach.

We present the computational complexity for the 2-D linear diffusion operator as an illustration. For the exact number of computations, we will have to isolate the interior and boundary cells to compute their individual number of operations and sum them up. Instead, since we are interested only in the order of magnitude of the computational complexity, we focus on computing the complexity for 1 cell and multiply that by the total number of cells, $N$. Consider the Equation 30 re-written here

$$\sum_f \left(\frac{\Gamma_f A_f}{\delta_f}\right)\left[\underbrace{(u_{k(f)} - u_0)}_{\text{I}} - \underbrace{\left[\frac{u_a - u_b}{|\vec{t}_f|}\right]\hat{t}_f.\vec{l}}_{\text{II}}\right]$$

The above term represents the discretization of the linear diffusion operator that has to be computed for every cell, where for each cell, the summation is over all the faces of the cell. Let $N$ denote the total number of cells, $N_f$ denote the total number of faces, $n_f$ denote the average number of faces per cell, and $n_c$ denote the average number of cells shared by a node.

The term $\left(\frac{\Gamma_f A_f}{\delta_f}\right)$ requires 2 operations and can be pre-computed for all the faces taking a total of $2 \times N_f$ operations. The I term is computed for every face in a cell and its product with the term within parentheses takes a total of $2 \times n_f \times N$. Now focusing on the term II, the dot product $\hat{t}_f.\vec{l}$ takes 3 operations (in 2D) and can be pre-computed for all the faces, leading to a total of $3 \times N_f$. The terms $u_a$ and $u_b$ have to be interpolated from cell-center values and its complexity depends on number of cells shared by a given node. The computation of each of the terms $u_a, u_b$ require $2n_c - 1$ operations. Therefore computing their difference, along with the multiplication and division operations of term II takes a total of $4n_c + 1$ and adding up multiplication with the term within parenthesis makes it $4n_c + 2$ per face and hence a total of $(4n_c + 2) \times n_f \times N$.

In total, the number of operations required to compute the linear diffusion operator is roughly $6N_f + [(4n_c + 4)n_f]N$, where the $6N_f$ arises from pre-computation of terms for every face and $[(4n_c + 4)n_f]N$ arises from computation of terms I and II. The number of cells shared by a node, $n_c$ and the number of faces per cell, $n_f$ are dependent on the mesh. Factoring this as a constant, and considering that the number of faces and cells in a mesh are of the same order of magnitude, the overall complexity is in $\mathcal{O}(N)$. Therefore, the cost of computing matrix $\mathbf{A}$ is linear with respect to the number of cells, $N$. Note that when $\mathbf{A}$ has parameter dependence, a unique $\mathbf{A}$ has to be computed for each snapshot and hence for $M$ snapshots, the overall offline cost for computing the linear operators is $\mathcal{O}(MN)$. This is considered as the price paid for the lack of access to the source code of the FOM, which is otherwise essential in the construction of a ROM, and is an offline cost used for model development. Finally, the RHS vector $\mathbf{f}$ needs to be constructed by applying the snapshot vectors to $\mathbf{A}$. Since $\mathbf{A}$ is known to have a highly sparse structure [50], the matrix-vector product would cost $\mathcal{O}(N)$ per snapshot leading to a total of $\mathcal{O}(MN)$. Therefore, the total cost for the offline construction of the matrices for the proposed methodology is still $\mathcal{O}(MN)$.

## IV. Results

The methodology described in this paper is demonstrated initially on a non-linear, parametric canonical PDE followed by the Euler equations simulating the flow past an airfoil.

### A. Canonical PDE

The canonical PDE is the same as that used in [37, 51] and is given below

$$-\nabla^2 u(x, y) + s(u(x, y); \mu) = 100sin(2\pi x)sin(2\pi y) \tag{39}$$

$$s(u; \mu) = \frac{\mu_1}{\mu_2}(e^{u\mu_2} - 1) \tag{40}$$

where, the spatial variables $(x, y) \in (0, 1)^2$ and the parameter $(\mu_1, \mu_2) \in [0.01,\ 2.0]^2$. The PDE is solved with Dirichlet boundary conditions of $u = 0$ along the boundaries. The Matlab PDE Toolbox [52] is the black-box code used in this test case to obtain the snapshots. The domain is discretized with an unstructured mesh comprising 1024 triangular cells ($N = 1024$). The computational grid and a sample snapshot solution of the PDE is shown in Figure 3. The following transformation is done to the PDE $[y_1,\ y_2] \to [u,\ \frac{\mu_1}{\mu_2}\,(e^{\mu_2 u} - 1)]$ which leads to the modified equation in terms of the observables

$$\begin{bmatrix} -\nabla^2 & 1 \end{bmatrix} \begin{bmatrix} y_1 \\ y_2 \end{bmatrix} = 100sin(2\pi x)sin(2\pi y) \tag{41}$$

where, $\nabla^2$ represents the laplacian. Upon discretization, the above equation transforms to

$$\begin{bmatrix} -\mathbf{L} & \mathbf{I} \end{bmatrix} \begin{bmatrix} \mathbf{y}_1 \\ \mathbf{y}_2 \end{bmatrix} + \mathbf{b}_a = 100sin(2\pi \mathbf{x})sin(2\pi \mathbf{y}) \tag{42}$$

where $\mathbf{L} \in \mathbb{R}^{N\times N}$ is the discrete laplacian and $\mathbf{I}$ is the identity matrix.

$$\underbrace{\begin{bmatrix} -\mathbf{L} & \mathbf{I} \end{bmatrix}}_{\mathbf{A}} \underbrace{\begin{bmatrix} \mathbf{y}_1 \\ \mathbf{y}_2 \end{bmatrix}}_{\mathbf{y}} = \underbrace{-\mathbf{b}_a + 100sin(2\pi \mathbf{x})sin(2\pi \mathbf{y})}_{\mathbf{f}} \tag{43}$$

which reduces to the $\mathbf{Ay} = \mathbf{f}$ form that we are interested in. The ROM is developed from this point using the same approach discussed in section II B.

In this test case, the observable $y_2$ is chosen to be $\frac{\mu_1}{\mu_2}(e^{\mu_2 u} - 1)$ and thereby lumping the parameter dependence along with the observable. This makes the linear operator independent of parameters and hence can be pre-computed. However, with a different choice of $y_2$ such as $e^{\mu_2 u}$, the linear operator, $\mathbf{A}$ becomes dependent on parameters, and therefore cannot be pre-computed. In such a case, a unique matrix is constructed for each snapshot, similar to the RHS, $\mathbf{f}$, and the linear operator for new parameter instances can be interpolated element-wise using methods in [46] or [47].

The constraint for the ROM for this case is given as follows

$$h(\mathbf{y}) = \mathbf{y}_2 - \frac{\mu_1}{\mu_2}(e^{\mu_2 \mathbf{y}_1} - 1) = 0 \tag{44}$$

The non-linear constraint equation above enforces the relationship between the two observables which should hold true. and is efficiently computed using the DEIM described in Appendix B.

In order to extract the POD basis, 20 snapshot locations were computed that vary the parameters $(\mu_1, \mu_2)$ using a Latin Hypercube Design [53]. The singular values of the resulting snapshot

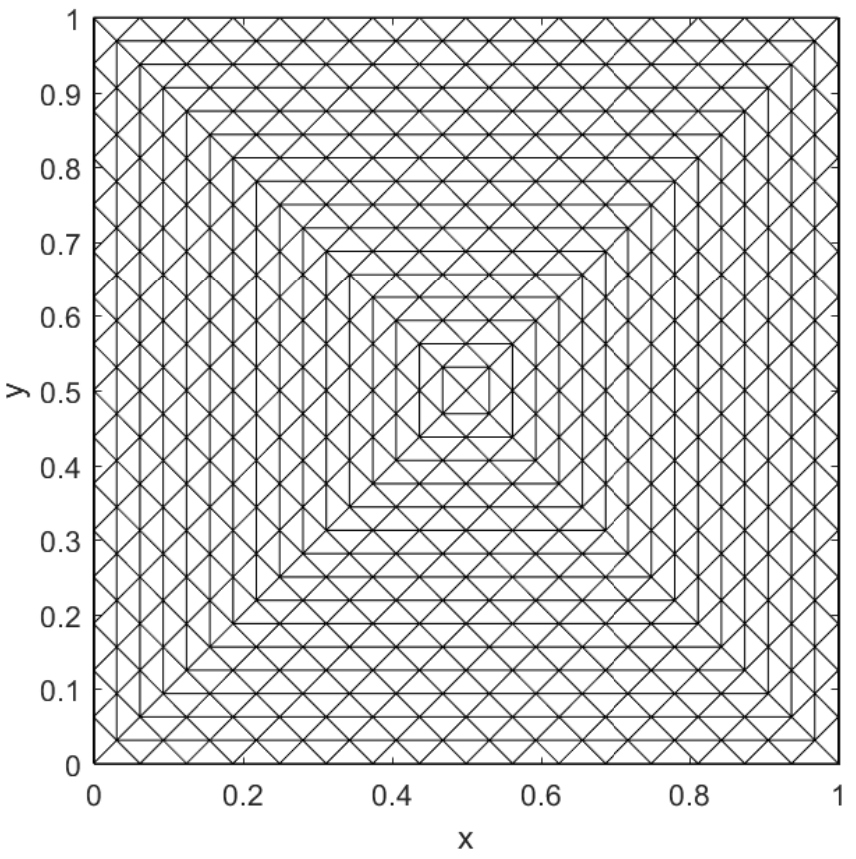

(a) Computational domain and mesh. 1024 triangular cells, 545 nodes and 1568 faces

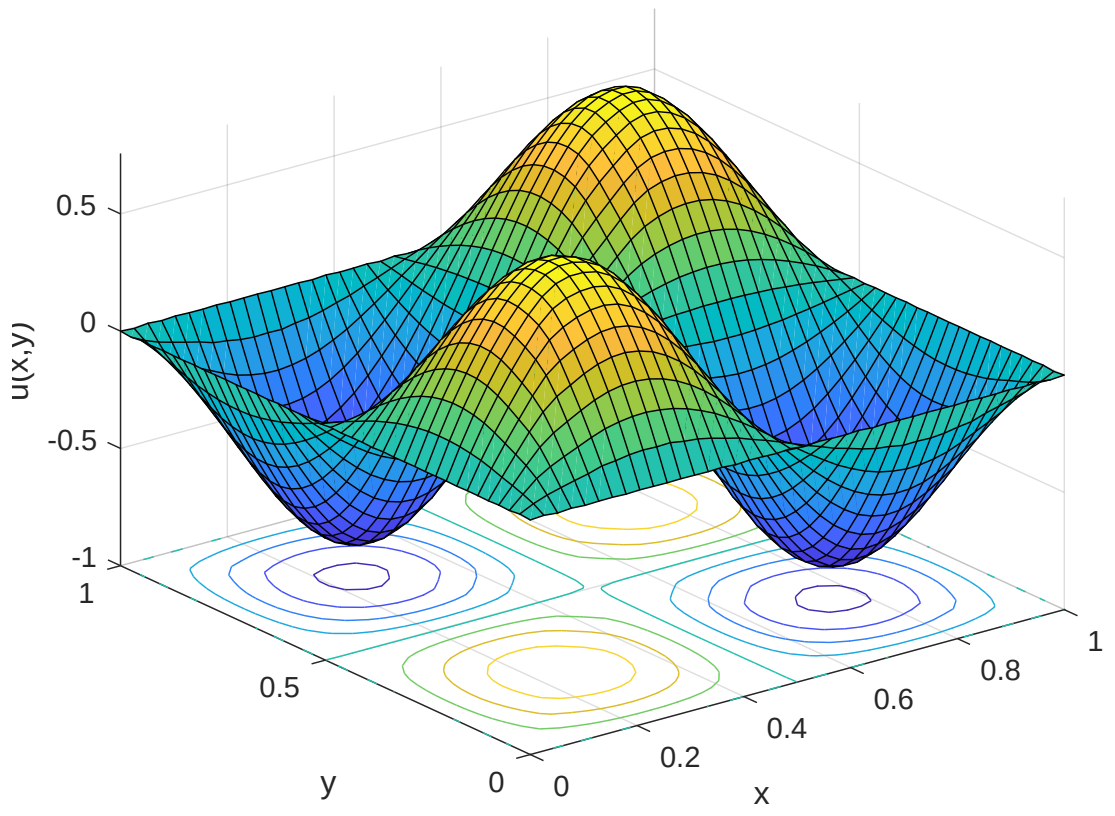

(b) Solution of non-linear test case at $\mu = (0.3, 9)$

**Figure 3 Computational domain and sample solution for the test cases**

matrix were not truncated to retain maximum accuracy. This results in a reduced system that is $20 \times 20$ which is still significantly smaller than the original system which is $1024 \times 1024$.

The comparison of the ROM and FOM solutions for a parameters outside of those used in the snapshots is shown in Figure 4. The linear operator matrix was pre-computed since it is independent of parameters while the RHS vector was piece-wise interpolated with a bi-variate $2^{nd}$ order lagrange polynomial. The relative error (R.E.) is computed as the relative value of the 2-norm error between the full and reduced order model solutions (Equation 45).

$$R.E. = \frac{\|FOM - ROM\|_2}{\|FOM\|_2} \tag{45}$$

The ROM results agree with a relative error in $\mathcal{O}(1)\%$ (see Table 1) which is an error measure based on the prediction through the entire computational domain, verifying that the proposed methodology is capable of accomplishing high levels of accuracy. The two main sources of error in this test case are (i) the approximation of the linear operator and (ii) the DEIM interpolation of the RHS vector. Since the Matlab PDE toolbox used to obtain the snapshots is a black-box it is not possible to quantify the error due to the linear operator approximation, and hence the overall error of approximation. However, given that the $R.E.$ for this test case was consistently in the $\mathcal{O}(1)\%$ which is a typically accepted range of accuracy of surrogate models used for engineering design optimization, it is concluded that the present method does approximately satisfy the governing equations at the ROM level, thereby serving as a physics-based surrogate model that is computationally cheap. As a next step, the method is applied towards approximating the flow past an airfoil, governed by the compressible Euler equations.

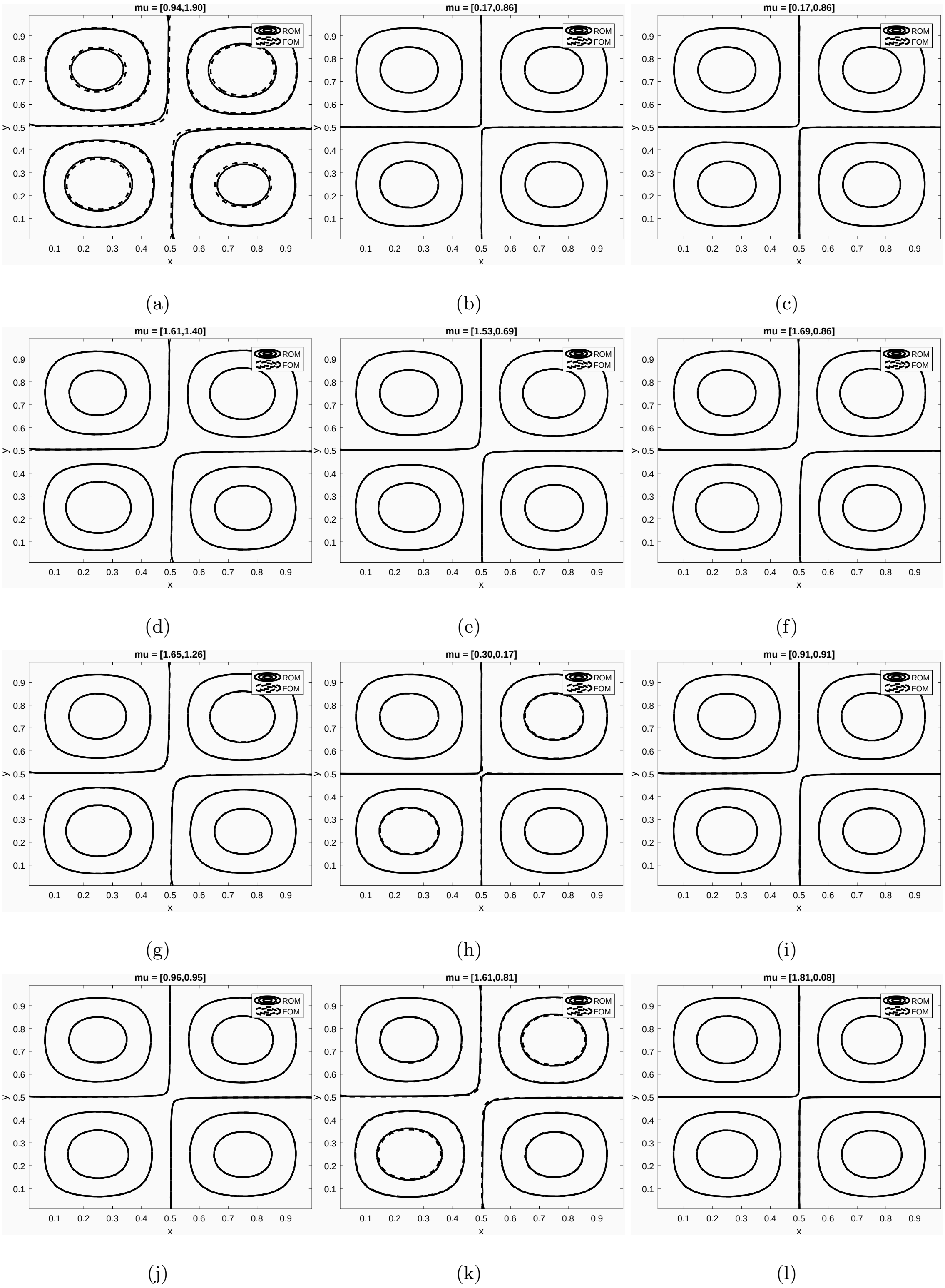

**Figure 4 Comparison of ROM (solid lines) and FOM (dashed lines) for the canonical PDE validation cases.**

**Table 1 Relative error of the ROM validation points for the Canonical PDE test case**

| Validation Case | $\mu_1$ | $\mu_2$ | R.E. | Validation Case | $\mu_1$ | $\mu_2$ | R.E. |
|---|---|---|---|---|---|---|---|
| 1 | 0.94 | 1.90 | 4.23 % | 7 | 1.65 | 1.26 | 0.52 % |
| 2 | 0.45 | 0.54 | 0.104 % | 8 | 0.30 | 0.17 | 1.08 % |
| 3 | 0.70 | 0.86 | 0.35 % | 9 | 0.91 | 0.91 | 0.125 % |
| 4 | 1.61 | 1.40 | 0.27 % | 10 | 0.96 | 0.95 | 0.17 % |
| 5 | 1.53 | 0.69 | 0.177 % | 11 | 1.61 | 0.81 | 2.05 % |
| 6 | 1.69 | 0.86 | 1.74 % | 12 | 1.81 | 0.08 | 0.11 % |

### B. Compressible Inviscid Flow past an Airfoil

The method is now implemented on the Euler equations governing the 2D, compressible, inviscid flow past an airfoil. The parameters varied are the Mach number ($\mathbb{M}$) and the angle of attack ($\alpha$), and the flow snapshots are generated by solving the PDEs in a commercial CFD solver, STARCCM+ [54]. The computational grid used by the CFD solver is exported in the CGNS [55] format and used to approximate the linear operator matrix in the present method via the finite volume method. The 2D compressible version of the Euler equations are given below in conservation form

$$\nabla_x \mathbf{F} + \nabla_y \mathbf{G} = 0 \tag{46}$$

where

$$\mathbf{F} = \begin{bmatrix} \rho u \\ \rho u^2 + p \\ \rho uv \\ \rho uH \end{bmatrix}, \ \mathbf{G} = \begin{bmatrix} \rho v \\ \rho uv \\ \rho v^2 + p \\ \rho vH \end{bmatrix}$$

$$H = E + \frac{p}{\rho}$$

$$\rho E = \frac{1}{2}\rho(u^2 + v^2) + \frac{p}{\gamma - 1}$$

and $\nabla_x$ and $\nabla_y$ are the $x$ and $y$ components of the gradient $\nabla$ respectively. The following transformation is performed

$$[\rho u, \rho v, \rho uv, p, \rho u^2, \rho v^2, \rho uH, \rho vH, \gamma]^T \to [y_1, y_2, y_3, y_4, y_5, y_6, y_7, y_8, y_9]^T$$

from the state variables to observables, leading to the transformed equation below. Note that in flows involving a calorically perfect gas, the $\gamma$ is a constant and hence its snapshots are not required to be collected. However, in the present study, $\gamma$ snapshots are collected for the sake of completeness and is used in the calculation of the energy term, $E$ in the non-linear constraints as shown in Equation 49.

$$\begin{bmatrix} \nabla_x & \nabla_y & & & & & & \\ & & \nabla_y & \nabla_x & \nabla_x & & & \\ & & \nabla_x & \nabla_y & & \nabla_y & & \\ & & & & & & \nabla_x & \nabla_y \end{bmatrix} \begin{bmatrix} y_1 \\ y_2 \\ y_3 \\ y_4 \\ y_5 \\ y_6 \\ y_7 \\ y_8 \end{bmatrix} = \mathbf{0} \tag{47}$$

The above equation upon discretization leads to

$$
\underbrace{\begin{bmatrix} \mathbf{G}_x & \mathbf{G}_y & & & & & & \\ & & \mathbf{G}_y & \mathbf{G}_x & \mathbf{G}_x & & & \\ & & \mathbf{G}_x & \mathbf{G}_y & & \mathbf{G}_y & & \\ & & & & & & \mathbf{G}_x & \mathbf{G}_y \end{bmatrix}}_{\mathbf{A}} \begin{bmatrix} \mathbf{y}_1 \\ \mathbf{y}_2 \\ \mathbf{y}_3 \\ \mathbf{y}_4 \\ \mathbf{y}_5 \\ \mathbf{y}_6 \\ \mathbf{y}_7 \\ \mathbf{y}_8 \end{bmatrix} = - \underbrace{\begin{bmatrix} \mathbf{b}_{a1} \\ \mathbf{b}_{a2} \\ \mathbf{b}_{a3} \\ \mathbf{b}_{a4} \\ \mathbf{b}_{a5} \\ \mathbf{b}_{a6} \\ \mathbf{b}_{a7} \\ \mathbf{b}_{a8} \end{bmatrix}}_{\mathbf{f}} \tag{48}
$$

where, $\mathbf{G}_x$ and $\mathbf{G}_y$ represents the discrete version gradient operators $\nabla_x$ and $\nabla_y$ respectively. With the FOM reduced to the $\mathbf{Ay} = \mathbf{f}$ form, the ROM is constructed as explained in section II B. The above set of equations are closed using non-linear constraints given by Equation 49. Notice that the constraints express the relationship between the first 4 observables ($y_1$ through $y_4$) and the last 4 observables ($y_5$ through $y_8$). It should be noted that all the observables that are in excess of the number of equations can be expressed as some function of the rest.

$$
\begin{aligned}
h_1 &= y_5 - \frac{y_1 y_3}{y_2} = 0 \\
h_2 &= y_6 - \frac{y_2 y_3}{y_1} = 0 \\
h_3 &= y_7 - y_1 \left( E + \frac{y_4 y_3}{y_1 y_2} \right) = 0 \\
h_4 &= y_8 - y_2 \left( E + \frac{y_4 y_3}{y_1 y_2} \right) = 0
\end{aligned} \tag{49}
$$

where $E = \frac{1}{2}\left(\frac{y_3 y_5}{y_1 y_2} + \frac{y_3 y_6}{y_1 y_2}\right) + \frac{y_3 y_4}{y_1 y_2 (y_9 - 1)}$. As mentioned before, the constraints are evaluated via the DEIM, explained in Appendix B. The Sequential Least Squares Quadratic Programming (SLSQP) [56] method available as part of the Matlab Non-Linear Optimization Toolbox is used to solve the resulting non-linear constrained optimization problem and the results are discussed as follows.

The NACA0012 is a symmetric airfoil as shown in Figure 5 and an unstructured mesh was generated to discretize a circular domain that is 150 chord lengths around it with a total of 11,265 triangular cells, 5776 vertices, and the near-field mesh is shown in Figure 5a. The parameters for this test case are the Mach number, $\mathbb{M}$ and angle of attack $\alpha$ and their ranges were set as $\mathbb{M} \in [0.3, 0.6]$ and $\alpha \in [0, 3]$ *deg.* 45 snapshots were generated using a Latin Hypercube Design out of which 40 were used to build the model while 5 were used to validate the model, see Figure 5b. The singular values that capture upto 99.99% of the variation of the observables were retained leading to a reduced order system of size $k = 148$.

**Table 2 Relative error of the ROM validation points for the NACA0012 test case**

| Validation Case | $\mathbb{M}$ | $\alpha$ [*deg.*] | Rel. Error |
|---|---|---|---|
| 1 | 0.51 | 1.77 | 1.008 % |
| 2 | 0.477 | 0.82 | 0.128 % |
| 3 | 0.32 | 1.36 | 0.69 % |
| 4 | 0.44 | 2.18 | 0.15 % |
| 5 | 0.46 | 0 | 0.194 % |

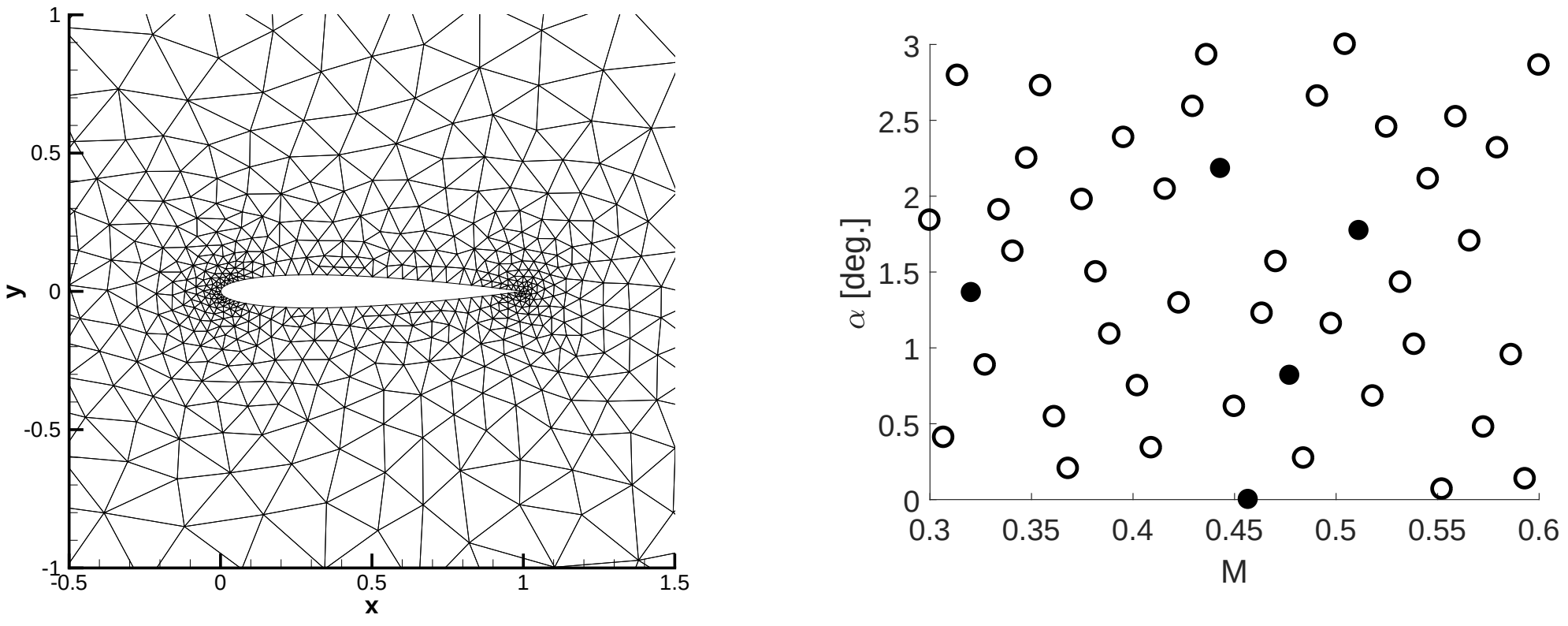

(a) Near field mesh for the NACA0012 test case

(b) Snapshot locations. Empty circles: Training points, Filled circles: Validation points

**Figure 5 Finite volume mesh and snapshot locations for the NACA0012 test case**

The relative error is calculated for the 5 validation points in Figure 5b and are summarized in Table 2. Note that these errors are calculated based on all 8 observables concatenated as a single vector. Overall, it is observed that the relative error of prediction is $\mathcal{O}(1)$ % similar to the canonical PDE test case.

The comparison of the ROM predictions against the actual FOM solutions are shown in Figure 6

in terms of the pressure and mach number contours. The plots show that the ROM prediction compares very well with the FOM results in terms of capturing the non-linear flow field. The accuracy of the ROM serves as a proof again that the current formulation satisfies governing equations at the ROM level, thereby serving as a physics-based surrogate model. Additionally, the ROM evaluates at a wall clock time in $\mathcal{O}(1)$ *sec* giving $\mathcal{O}(100\times)$ speedup compared to the FOM on a desktop computer, offering the suitability to be used in a design optimization framework.

### C. Discussion

The training data is an important aspect of any surrogate model development and for model order reduction, this relates to the location and number of snapshots used. For the present study, 20 and 40 snapshots were used respectively for the canonical PDE and the NACA0012 test cases. Also, the snapshots are chosen via a 'space-filling' design that maximizes the minimum distance between their locations in the design space. Based on the model validation conducted, the number and locations of the snapshots are justified for the present study, but the authors do acknowledge that a sequential approach that adaptively samples snapshots and updates the model would be necessary to scale the method to large dimensional problems.

A key step in parametric model reduction is the interpolation of the ROM system matrices for new parameter instances and the accuracy of the interpolation is critical to the overall accuracy of the method. In this work, as mentioned earlier, a piece-wise interpolation in the original space of the reduced system matrices was performed using polynomials in the Lagrange form. While in the present study, it was observed that a $2^{nd}$ order interpolation was sufficient to get ROM accuracy in $\mathcal{O}(1)\%$, this might not always be enough for the same problem for instance when new parameters are added or the existing parameter ranges are expanded. Therefore, the authors acknowledge that the choosing the most suitable interpolation method for a given problem is in principle a part of the method itself. However, that aspect is not explored in the current work to focus purely on demonstrating the methodology.

The flow physics chosen for the current test problem was restricted to 2D, inviscid flow and the design variable ranges were chosen such that the FOM is restricted to a non-linear but shock-free and small angles of attack flow regime, so that a smaller computational mesh suitable for preliminary

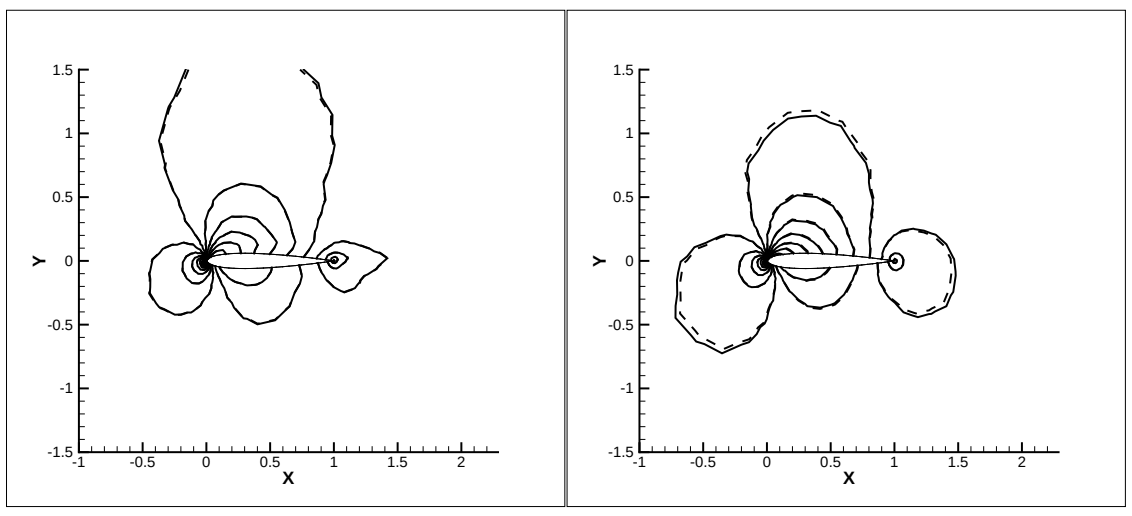

(a) $\mathbb{M} = 0.51$ $\alpha = 1.77$ *deg.*

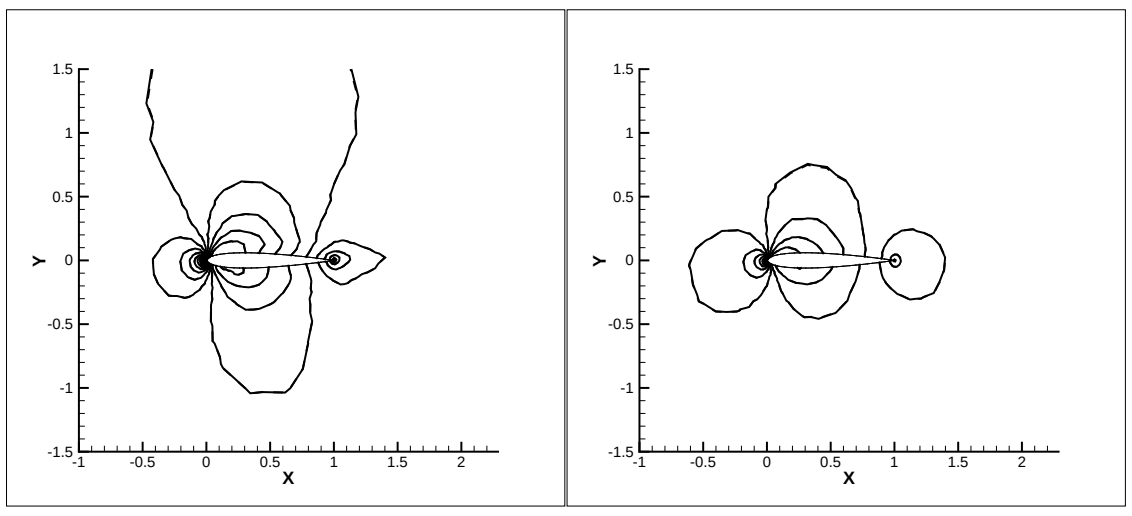

(b) $\mathbb{M} = 0.477$ $\alpha = 0.82$ *deg.*

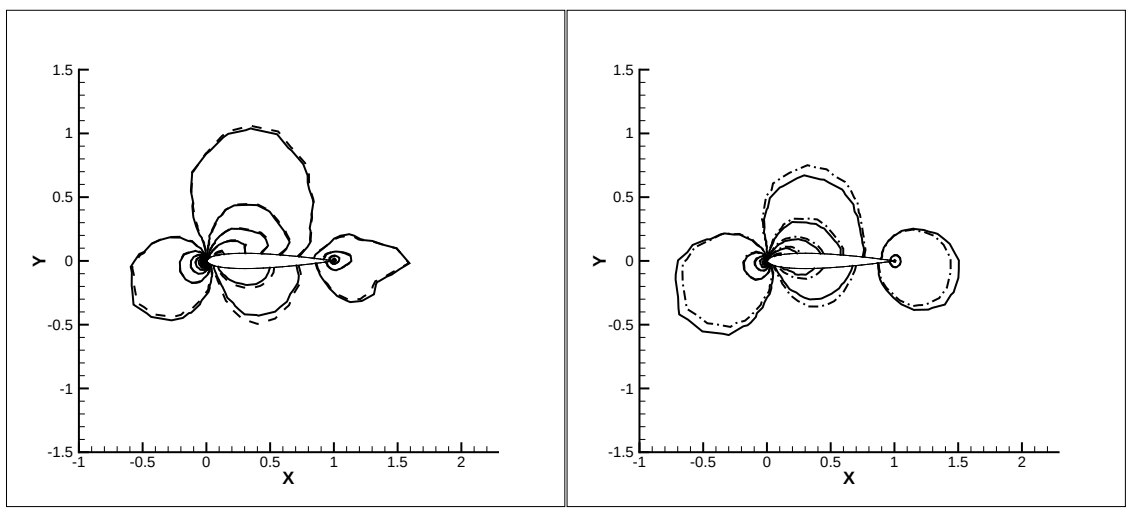

(c) $\mathbb{M} = 0.32$ $\alpha = 1.36$ *deg.*

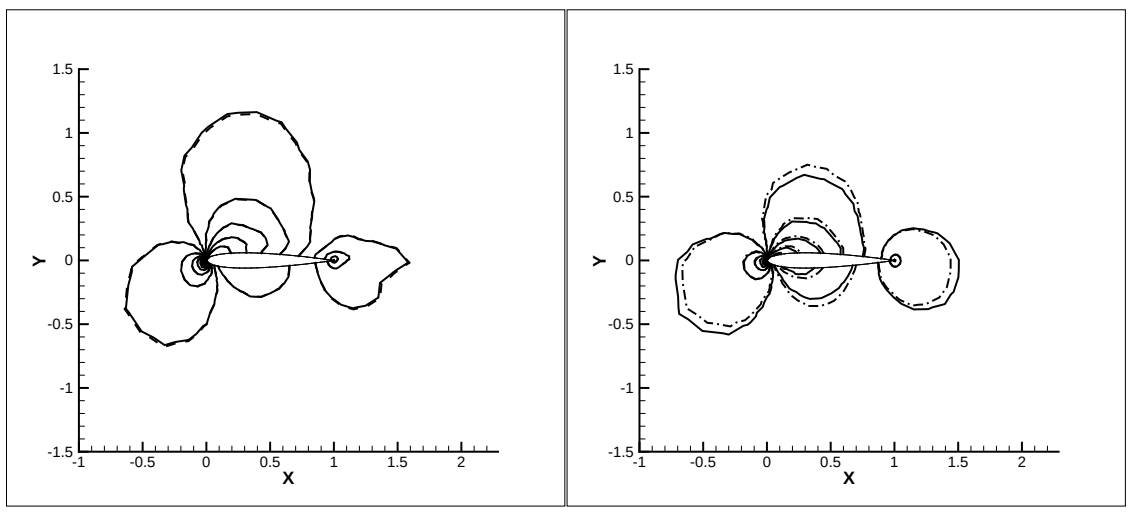

(d) $\mathbb{M} = 0.44$ $\alpha = 2.18$ *deg.*

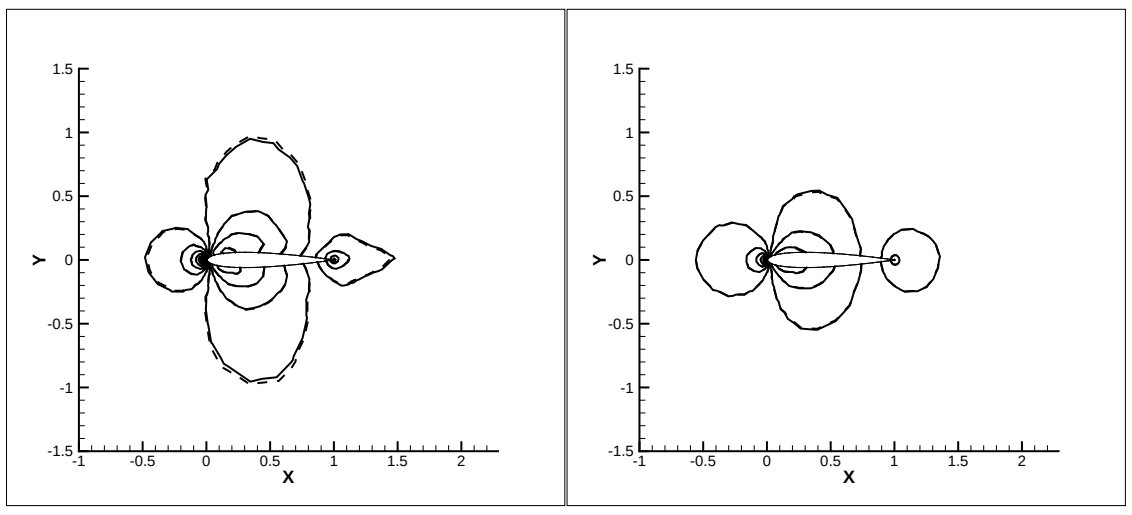

(e) $\mathbb{M} = 0.46$ $\alpha = 0$ *deg.*

**Figure 6 Comparison of ROM (dashed lines) and FOM (solid lines) for the NACA0012 validation cases. In each figure, left=Mach contours, right=Pressure contours**

evaluation of the method was possible. However, the method equally extends to 3D viscous flows and flows in the transonic and supersonic regimes without modification. Additionally, the design parameters chosen for this test case are boundary parameters and hence are lumped alongside the RHS term, resulting in a matrix, $\mathbf{A}$ that is independent of the parameters and hence its projection can be pre-computed. However, with the addition of other types of parameters such as geometry (or shape) parameters that lead to a unique $\mathbf{A}$ for each snapshot, the proposed method still applies without modification.

## V. Conclusion

This paper proposes a method that enables projection-based model reduction with black-box full order models governed by steady, non-linear parametric PDEs that are widely used in aerospace design optimization. It particularly addresses the situation where there is knowledge but no access to the mathematical model such as in the case of commercial CFD codes. It does so by retrieving the discretized form of the governing equations from system outputs and domain knowledge. The method first transforms the state variables and their non-linear terms into a set of linear observables resulting in a linear equation. The method then approximates the linear differential operator matrix through a direct finite volume discretization, which requires the computational grid and knowledge of the governing equations. Even though a linear transformation is performed on the governing equations, a set of non-linear constraints are required to ensure the equations of state are satisfied and their computation is efficiently handled via the DEIM.

The method relies on generating a database of ROMs corresponding to several pre-determined set of parameter combinations which are later interpolated for new parameter instances. An accurate interpolation method for the ROMs is key to the success of this approach, but the current paper focuses only on demonstrating the overall methodology. The methodology was successfully tested on a 2D canonical PDE with 2 parameters and exponential non-linearity and the 2D compressible Euler equations with 2 parameters and polynomial-type non-linearities. In either case the ROM accuracy were observed in $\mathcal{O}(1)$ % while the computational speedups were in $\mathcal{O}(100\times)$ making it suitable for use in design optimization and uncertainty quantification applications.

The proposed approach comes with an offline cost to construct the linear operators in $\mathcal{O}(N)$ per snapshot. This cost ignores a pre-processing step in the finite volume discretization where geometric parameters such as the surface normals, inter-cell distances and interpolants are to be computed for each cell in the grid, since certain grid generators or CFD packages already provide this information. Due to the lack of access to the source code of the FOM, the authors consider such a cost inevitable in enabling projection-based model reduction with black-box FOMs. However, the method enables the development of a physics-based surrogate model that satisfies a reduced form of the actual governing equations which is essential in the design of next generation aerospace concepts. Having established feasibility of the proposed method, the logical next steps are to expand the scope of the parameters (transonic/supersonic regimes) and physics (3D/viscous flows) and demonstrate the applicability of the method on suitable design optimization problems.

**Acknowledgements**

Professor Yingjie Liu acknowledges the support in part by NSF grants DMS-1522585 and DMS-1622453.

## References

[1] Forrester, A. I. J. and Keane, A. J., “Recent advances in surrogate-based optimization,” *Progress in Aerospace Sciences*, Vol. 45, No. 1-3, 2009, pp. 50–79.

[2] Queipo, N. V., Haftka, R. T., Shyy, W., Goel, T., Vaidyanathan, R., and Kevin Tucker, P., “Surrogate-based analysis and optimization,” *Progress in Aerospace Sciences*, Vol. 41, No. 1, 2005, pp. 1–28.

[3] Viana, F. A. C., Haftka, R. T., and Watson, L. T., “Efficient global optimization algorithm assisted by multiple surrogate techniques,” *Journal of Global Optimization*, Vol. 56, No. 2, 2013, pp. 669–689.

[4] Andy J Keane, P. B. N., *Computational Approaches for Aerospace Design: The Pursuit of Excellence*, WILEY, 2005.

[5] Bui-Thanh, T., Damodaran, M., and Willcox, K. E., “Aerodynamic Data Reconstruction and Inverse Design Using Proper Orthogonal Decomposition,” *AIAA Journal*, Vol. 42, No. 8, 2004, pp. 1505–1516.

[6] Legresley, P. a. and Alonso, J. J., “Airfoil Design Optimization Using Reduced Order Models Based on Proper Orthogonal Decomposition FLUIDS,” *AIAA Fluids 2000 Conference and Exhibit*, 2000.

[7] Bui-Thanh, T., Willcox, K., and Ghattas, O., “Parametric Reduced-Order Models for Probabilistic Analysis of Unsteady Aerodynamic Applications,” *AIAA Journal*, Vol. 46, No. 10, 2008, pp. 2520–2529.

[8] Haddadpour, H., Behbahani-Nejad, M., and Firooz-Abadi, R. D., “Reduced Order Aerodynamic Model for Aeroelastic Analysis of Complex Configurations in Incompressible Flow,” *Journal of Aircraft*, Vol. 44, No. 3, 2007, pp. 1015–1019.

[9] Lieu, T. and Farhat, C., “Adaptation of Aeroelastic Reduced-Order Models and Application to an F-16 Configuration,” *AIAA Journal*, Vol. 45, No. 6, 2007, pp. 1244–1257.

[10] Brouwer, K. R., Crowell, A. R., and Mcnamara, J. J., “Rapid Prediction of Unsteady Aeroelastic Loads in Shock-Dominated Flows,” , No. January, 2015, pp. 1–20.

[11] Dowell, E. and Hall, K., “Reduced Order Models in Unsteady Aerodynamic Models, Aeroelasticity and Molecular Dynamics,” *ICAS - 26th Congress of International Council of the Aeronautical Sciences 2008*, 2006, pp. 1–13.

[12] Glaz, B., Liu, L., and Friedmann, P. P., “Reduced-Order Nonlinear Unsteady Aerodynamic Modeling Using a Surrogate-Based Recurrence Framework,” *AIAA Journal*, Vol. 48, No. 10, 2010, pp. 2418–2429.

[13] Karpel, M., “Reduced-Order Models for Integrated Aeroservoelastic Optimization,” *Journal of Aircraft*, Vol. 36, No. 1, 1999.

[14] Lewin, G. C. and Haj-Hariri, H., "Reduced-Order Modeling of a Heaving Airfoil," *AIAA Journal*, Vol. 43, No. 2, 2005, pp. 270–283.

[15] Skujins, T. and Cesnik, C. E. S., "Reduced-Order Modeling of Unsteady Aerodynamics Across Multiple Mach Regimes," *Journal of Aircraft*, Vol. 51, No. 6, 2014, pp. 1681–1704.

[16] Thomas, J. P., Dowell, E. H., and Hall, K. C., "Three-Dimensional Transonic Aeroelasticity Using Proper Orthogonal Decomposition-Based Reduced-Order Models," *Journal of Aircraft*, Vol. 40, No. 3, 2003, pp. 544–551.

[17] Crowell, A. R., Mcnamara, J. J., Kecskemety, K. M., and Goerig, T. W., "A Reduced Order Aerothermodynamic Modeling Framework for Hypersonic Aerothermoelasticity," *Direct*, Vol. 2969, No. April, 2010, pp. 22.

[18] Falkiewicz, N. J., S. Cesnik, C. E., Crowell, A. R., and McNamara, J. J., "Reduced-Order Aerothermoelastic Framework for Hypersonic Vehicle Control Simulation," *AIAA Journal*, Vol. 49, No. 8, 2011, pp. 1625–1646.

[19] Buffoni, M. and Willcox, K., "Projection-based model reduction for reacting flows," *40th Fluid Dynamics Conference and Exhibit*, No. July, 2010, pp. 1–14.

[20] Chen, P. C., Liu, D. D., and Chang, K. T., "TPS Design by a POD / RSM-Based Approach," *Aerospace Engineering*, , No. January, 2006, pp. 1–24.

[21] Qian, J., Wang, Y., Song, H., Pant, K., Peabody, H., Ku, J., and Butler, C. D., "Projection-Based Reduced-Order Modeling for Spacecraft Thermal Analysis," *Journal of Spacecraft and Rockets*, Vol. 52, No. 3, 2015, pp. 978–989.

[22] Epureanu, B. I., Dowell, E. H., and Hall, K. C., "Mach Number Influence on Reduced-Order Models of Inviscid Potential Flows in Turbomachinery," *Journal of Fluids Engineering*, Vol. 124, No. 4, 2002, pp. 977–987.

[23] Epureanu, B. I., "A parametric analysis of reduced order models of viscous flows in turbomachinery," *Journal of Fluids and Structures*, Vol. 17, No. 7, 2003, pp. 971–982.

[24] Florea, R., Hall, K. C., and Cizmas, P. G. a., "Reduced-order modeling of unsteady viscous flow in a compressor cascade," *AIAA Journal*, Vol. 36, No. 6, 1998, pp. 1039–1048.

[25] Fossati, M., "Evaluation of Aerodynamic Loads via Reduced-Order Methodology," *AIAA Journal*, Vol. 53, No. 8, 2015, pp. 1–17.

[26] Gennaretti, M. and Greco, L., "Time-Dependent Coefficient Reduced-Order Model for Unsteady Aerodynamics of Proprotors," *Journal of Aircraft*, Vol. 42, No. 1, 2005, pp. 138–147.

[27] Gennaretti, M. and Muro, D., "Multiblade Reduced-Order Aerodynamics for State-Space Aeroelastic Modeling of Rotors," *Journal of Aircraft*, Vol. 49, No. 2, 2012, pp. 495–502.

[28] Lucia, D. J., King, P. I., and Beran, P. S., "Domain Decomposition for Rediced-Order Modeling of a Flow with Moving Shocks," *AIAA Journal*, Vol. 40, No. 11, 2002, pp. 2360–2363.

[29] Lucia, D. J. and Beran, P. S., "Projection methods for reduced order models of compressible flows," *Journal of Computational Physics*, Vol. 188, No. 1, 2003, pp. 252–280.

[30] Lucia, D. J., King, P. I., and Beran, P. S., "Reduced order modeling of a two-dimensional flow with moving shocks," *Computers and Fluids*, Vol. 32, No. 7, 2003, pp. 917–938.

[31] Lucia, D. J., Beran, P. S., and Silva, W. a., "Reduced-order modeling: New approaches for computational physics," *Progress in Aerospace Sciences*, Vol. 40, No. 1-2, 2004, pp. 51–117.

[32] Antoulas, A. C. and Sorensen, D. C., "Approximation of large-scale dynamical systems : An overview," Vol. 1892, 2001, pp. 1–22.

[33] Benner, P., Gugercin, S., and Willcox, K., "A Survey of Model Reduction Methods for Parametric Systems," *Max Planck Institute Magdeburg*, 2013.

[34] Holmes, Philip., Lumley, John L., Berkooz, Gahl and Rowley, C. W., *Turbulence, Coherent Structures, Dynamical Systems and Symmetry*, Vol. 36, 1998.

[35] Holmes, Philip., Lumley, John L., Berkooz, Gahl and Rowley, C. W., *Turbulence, Coherent Structures, Dynamical Systems and Symmetry*, Vol. 36, 1998.

[36] Chatterjee, A., "An introduction to the proper orthogonal decomposition," *Current Science*, Vol. 78, No. 7, 2000, pp. 808–817.

[37] Chaturantabut, S. and Sorensen, D. C., "Nonlinear Model Reduction via Discrete Empirical Interpolation," *SIAM Journal on Scientific Computing*, Vol. 32, No. 5, 2010, pp. 2737–2764.

[38] Xiao, D., Fang, F., Buchan, A. G., Pain, C. C., Navon, I. M., and Muggeridge, A., "Non-intrusive reduced order modelling of the Navierâ??Stokes equations," *Computer Methods in Applied Mechanics and Engineering*, Vol. 293, 2015, pp. 522–541.

[39] D. Xiao, F. Fang, C. P. and Hu, G., "Non-intrusive reduced-order modelling of the Navier–Stokes equations based on RBF interpolation," *INTERNATIONAL JOURNAL FOR NUMERICAL METHODS IN FLUIDS*, Vol. 79, 2015, pp. 580–595.

[40] Kutz, N., Proctor, J. L., and Brunton, S. L., "Koopman Theory for Partial Differential Equations," *arXiv:1607.07076v1*, 2016, pp. 1–21.

[41] Williams, M. O., Kevrekidis, I. G., and Rowley, C. W., "A Data-Driven Approximation of the Koopman Operator: Extending Dynamic Mode Decomposition," *Journal of Nonlinear Science*, Vol. 25, No. 6,

2015, pp. 1307–1346.

[42] Peherstorfer, B. and Willcox, K., "Dynamic data-driven reduced-order models," *Computer Methods in Applied Mechanics and Engineering*, Vol. 291, 2015, pp. 21–41.

[43] Kutz, J. N., Proctor, J. L., and Brunton, S. L., "GENERALIZING KOOPMAN THEORY TO ALLOW FOR INPUTS AND CONTROL," *arXiv:1602.07647v1*, 2016, pp. 1–21.

[44] Brunton, S. L., Brunton, B. W., Proctor, J. L., and Kutz, J. N., "Koopman invariant subspaces and finite linear representations of nonlinear dynamical systems for control," *PLoS ONE*, Vol. 11, No. 2, 2016, pp. 1–19.

[45] Versteeg,H.K., Malalasekera, W., *An Introduction to Computational Fluid Dynamics - The Finite Volume Method*, Prentice Hall, 2007.

[46] Joris Degroote, J. V. and Willcox, K., "Interpolation among reduced-order matrices to obtain parameterized models for design, optimization and probabilistic analysis," *International Journal for Numerical Methods in Fluids*, Vol. 63, 2010, pp. 207–230.

[47] Amsallem, D. and Farhat, C., "Interpolation Method for Adapting Reduced-Order Models and Application to Aeroelasticity," *AIAA Journal*, Vol. 46, No. 7, 2008, pp. 1803–1813.

[48] Jeffreys, H. and Jeffreys, B. S., *Methods of Mathematical Physics*, chap. Lagrange's Interpolation Formula., Cambridge University Press, 3rd ed., 1988, p. 260.

[49] Krady, M. M., "Extension Of Lagrange Interpolation," *INTERNATIONAL JOURNAL OF SCIENTIFIC & TECHNOLOGY RESEARCH*, 2015.

[50] Renganathan, S. A. and Mavris, D. N., "A Methodology for Projection-Based Model Reduction with Black-Box High-Fidelity Models," *17th AIAA Aviation Technology, Integration, and Operations Conference*, 2017, p. 4444.

[51] Patera, A. T., "Reduced Basis Approximation and A Posteriori Error Estimation for Parametrized Partial Differential Equations," *Foundations*, , No. January, 2007, pp. 251.

[52] "MATLAB PDE Toolbox Release R2017b, The MathWorks Inc." 2017.

[53] Myers, R. H., Montgomery, D. C., and Anderson-Cook, C. M., *Response surface methodology: process and product optimization using designed experiments*, John Wiley & Sons, 2016.

[54] "STARCCM+ url: http://mdx.plm.automation.siemens.com/star-ccm-plus," 2017.

[55] Poirier, D., Allmaras, S. R., McCarthy, D. R., Smith, M. F., and Enomoto, F. Y., "The CGNS system," *AIAA paper*, , No. 98-3007, 1998.

[56] Kraft, D. and Schnepper, K., "SLSQP A Nonlinear Programming Method with Quadratic Programming Subproblems," *DLR, Oberpfaffenhofen*, 1989.

### Appendix A: Koopman Theory

The Koopman Theory of PDEs [40] states that upon judicious choice of a set of scalar observables that are non-linear functions of the state variables, the original finite-dimensional non-linear dynamical system can be represented as an infinite-dimensional linear system. That is, given the following evolution of the state for a non-linear dynamical system

$$\mathbf{u}_{t+1} = \mathbf{F}(\mathbf{u}_t) \tag{A1}$$

The Koopman operator is defined as

$$\mathcal{K}g(\mathbf{u}_t) = g(\mathbf{F}(\mathbf{u}_t)) = g(\mathbf{u}_{t+1}) \tag{A2}$$

Therefore, an infinite dimensional but linear map, $\mathcal{K}$ can be found for a non-linear system. Further, it has been shown that a finite dimensional approximation,$K$ to the infinite dimensional $\mathcal{K}$ can be determined through a least-squares minimization problem that considers a dictionary of observables [41]. i.e. by collecting snapshots of an assumed set of observables, the operator $K$ can be determined as the best-fit operator between the observables at time $t$ and $t+1$ in the least squares sense. Additionally, [42] have shown that such a method converges to the actual governing equations (in terms of the linear and non-linear operators) with sufficient trajectory data.

This work draws from the Koopman theory and represents the underlying steady non-linear parametric system as a linear system in terms of the observables. In this regard, the present method assumes knowledge of the governing equations of the FOM, although no access to their discretized form in the source code is necessary. Since there is no trajectory data, the linear operator is directly approximated by discretizing the linear differential terms through a suitable discretization method such as the finite volume method.

### Appendix B: Discrete Empirical Interpolation Method (DEIM)

The Discrete Empirical Interpolation Method (DEIM) is briefly reviewed here and as an illustration one of the non-linear constraints used in the Airfoil test case is evaluated.

For a non-linear function $\mathbf{f}(\theta) \in \mathbb{R}^N$ the DEIM approximates $\mathbf{f}$ by projecting it onto a subspace spanned by $\{\mathbf{x}_1, ..., \mathbf{x}_q\} \subset \mathbb{R}^N$ as

$$\mathbf{f}(\theta) \approx \mathbf{X}c(\theta) \tag{B1}$$

where $\mathbf{X} = [\mathbf{x_1}, ..., \mathbf{x_q}] \in \mathbb{R}^{N \times q}, \quad q << N$ is determined via a POD of the snapshots of $\mathbf{f}$ and is assumed to be globally valid in the design space that bounds the design parameters $\theta$ and $\mathbf{c}(\theta) \in \mathbb{R}^q$ are the coefficients of the basis expansion. Then the approximation of $\mathbf{f}$ requires only the determination of $\mathbf{c}(\theta)$ which requires only $q$ equations. The DEIM gives a distinguished set of $q$ points from the over-determined system $\mathbf{f}(\theta) = \mathbf{X}\mathbf{c}(\theta)$. Given a permutation matrix $\mathbf{P}$ that would give $q$ such distinguished rows of a matrix when pre-multiplied, then the $q \times q$ system necessary to solve for the coefficients is given by

$$\mathbf{P}^T\mathbf{f}(\theta) = (\mathbf{P}^T\mathbf{X})\mathbf{c}(\theta) \tag{B2}$$

So the approximation of $\mathbf{f}(\theta)$ is then given by

$$\mathbf{f}(\theta) \approx \mathbf{X}(\mathbf{P}^T\mathbf{X})^{-1}\mathbf{P}^T\mathbf{f}(\theta) \tag{B3}$$

If the $q$ row-indices (that are extracted by pre-multiplying with $\mathbf{P}^T$) are represented by a vector, $\varrho$, then in the above equation, $\mathbf{P}^T\mathbf{f}(\theta)$ is equivalent to extracting the $\varrho$ rows of $\mathbf{f}$. Therefore the approximation of $\mathbf{f}(\theta)$ requires only $q$ computations which is efficient because $q << N$. Similarly, a non-linear function that depends on the state, $\mathbf{f}(\mathbf{u})$ can be approximated as

$$\mathbf{f}(\mathbf{u}) \approx \mathbf{X}(\mathbf{P}^T\mathbf{X})^{-1}\mathbf{P}^T\mathbf{f}(\mathbf{u}) \tag{B4}$$

Since $\mathbf{u} = \Phi_k^T\tilde{\mathbf{u}}$ and setting $\tilde{\mathbf{f}} = \Phi_k^T\mathbf{f}(\mathbf{u})$, $\tilde{\mathbf{f}}$ can be approximated as

$$\tilde{\mathbf{f}} = \Phi_k^T\mathbf{X}(\mathbf{P^T}\mathbf{X})^{-1}\mathbf{f}(\mathbf{P}^T\Phi_k\tilde{\mathbf{u}}) \tag{B5}$$

In the above equation, the term $\Phi_k^T\mathbf{X}(\mathbf{P^T}\mathbf{X})^{-1}$ is independent of the state and hence can be pre-computed and $\mathbf{P}^T\Phi_k$ is just extraction of the $\varrho$ rows of $\Phi_k$. Therefore using the DEIM, the non-linear can be expressed in terms of the reduced state, $\tilde{\mathbf{u}}$ and hence can be efficiently computed.

Now the DEIM is illustrated on evaluating the first constraint of Equation 49 which in discretized form is given below

$$\mathbf{h}_1 = \mathbf{y}_5 - \frac{\mathbf{y}_1\mathbf{y}_3}{\mathbf{y}_2} \tag{B6}$$

Let $\varrho_5$ be the vector containing the $q$ row-indices returned by DEIM via snapshots of the non-linear term $\mathbf{y}_5$ and $\phi_1$, $\phi_2$, $\phi_3$, $\phi_5$ be the projection matrix of $\mathbf{y}_1$, $\mathbf{y}_2$, $\mathbf{y}_3$ and $\mathbf{y}_5$ respectively. Then

$$\tilde{\mathbf{h}}_1 = \tilde{\mathbf{y}_5} - \phi_5^T \mathbf{X} \left[\mathbf{X}(\varrho_5, :)\right]^{-1} \left\{ \frac{\phi_1(\varrho_5, :)\tilde{\mathbf{y}_1} \quad \phi_3(\varrho_5, :)\tilde{\mathbf{y}_3}}{\phi_2(\varrho_5, :)\tilde{\mathbf{y}_2}} \right\} \tag{B7}$$

In the above equation, the term outside of the braces can be pre-computed. Additionally since $\mathbf{y}_5 = \frac{\mathbf{y}_1\mathbf{y}_3}{\mathbf{y}_2}$, $\mathbf{X} = \phi_5$ and hence the term reduces to $[\mathbf{X}(\varrho_5, :)]^{-1}$ which is $q \times q$ and hence can be cheaply computed. Therefore using the DEIM, the non-linear constraints are evaluated in terms of the reduced state variables which makes it computationally cheap.